\documentclass[11pt]{article}

\usepackage{amsmath,amsthm,amssymb,amsfonts,amscd,array}
\usepackage{url}

\newtheorem{defn}{Definition}
\newtheorem{prop}{Proposition}

\title{On the Smallest Cayley Table Fill That Uniquely Determines a Group}
\author{%
Andrey M. Chernev\thanks{Corresponding author. Email: chernevAM@mpei.ru}\\
National Research University ``Moscow Power Engineering Institute''
\and
Nikolay A. Chernev\thanks{Email: nachernev@gmail.com}\\
Bauman Moscow State Technical University
}
\date{}

\begin{document}

\maketitle

\begin{abstract}
The article investigates the minimal number of cells in a Cayley table that must be filled in order to uniquely determine a group. An exact result is obtained for groups of orders \(2\)--\(23\), and an upper bound is given for an arbitrary order.
\end{abstract}



\section{Introduction}

One of the efficient tools for analysing a finite algebraic structure is the Cayley table (multiplication table), for example, for detecting regular elements and idempotents, finding subsemigroups, etc.

It is interesting to ask how much information about the properties of an algebraic structure can be extracted from a partially filled Cayley table. In particular, in [1] the question was studied for which fillings of the table the proportion of associative groupoids (i.e.\ semigroups) is minimal or maximal. In [1], an example of a filling of the Cayley table was constructed which excludes associativity (i.e.\ the proportion of semigroups is \(0\)), and it was noted that a more interesting modified problem is the following: for which partial filling of the Cayley table is the proportion of semigroups minimal but not equal to \(0\).

The question can be reformulated as follows: suppose that we know that a partially filled Cayley table defines a semigroup. What is the minimum number of cells that must be filled so that this semigroup is determined uniquely (up to isomorphism)? This question was studied in [2], where an answer was obtained for semigroups of orders \(2\)--\(7\).

In the present article we study the analogous question for groups.

The obtained results are included in the Online Encyclopedia of Integer Sequences (OEIS), the largest database identifying and systematizing patterns in integers (for more details on OEIS see [3]).

\section{Preliminary information}

Let us recall some definitions used in this article.

\begin{defn}
A \emph{group} is a set endowed with a binary operation satisfying the following properties:
\begin{itemize}
\item associativity, i.e.\ for any elements \(a,b,c\),
\[
a(bc)=(ab)c;
\]
\item existence of a neutral element (identity): there exists an element \(e\) such that for every \(a\),
\[
ae=ea=a;
\]
\item existence of an inverse element: for every element \(a\) of the group there exists an element \(a^{-1}\) such that
\[
aa^{-1}=a^{-1}a=e.
\]
\end{itemize}
\end{defn}

\begin{defn}
The \emph{Cayley table} (multiplication table) of a group is the table composed of the results of applying the binary operation to the elements of the group; in the \(i\)-th row and the \(j\)-th column there is the product of the \(i\)-th and \(j\)-th elements of the group. In what follows, the entry of the Cayley table in the \(i\)-th row and the \(j\)-th column will be denoted by \(T_{ij}\).
\end{defn}

\begin{defn}
The \emph{order} of a finite group is the number of its elements.
\end{defn}

\begin{defn}
The \emph{order of an element} of a group (denoted by \(\operatorname{ord}(a)\)) is the minimum \(n\) for which
\[
a^n=e.
\]
It is well known that the order of an element divides the order of the group.
\end{defn}

\begin{defn}
A group of order \(n\) consisting only of the elements
\[
a,a^2,\ldots,a^{n-1},\qquad a^n=e
\]
is called \emph{cyclic}.
\end{defn}

\section{Methodology of the research}

The character of the problem of the minimal filling of the Cayley table changes in comparison with semigroups, where the problem is purely combinatorial.

First of all, in the case when the order of the group is coprime to \(\varphi(n)\) (where \(\varphi\) is the Euler function), there exists only one group of order \(n\) (this result was obtained in [4]). In this case, no filling of the Cayley table is required for the unique identification of the group, i.e.\ the minimum number of filled cells is equal to \(0\).

Even if there are several groups of a given order, their number is small, and for groups there is a constructive method for finding a minimal filling, namely: one formulates a property distinguishing the given group from the others and expresses this property with the help of cells of the Cayley table.

The property ``not to be commutative'' is expressed in the most economical way.

\begin{prop}
To characterize a group as nonabelian, it is sufficient to fill two cells of the Cayley table.
\end{prop}

\begin{proof}
Choose two noncommuting elements and put them in positions \(0\) and \(1\) in the enumeration of the elements of the group. Fill the Cayley table as follows:
\[
T_{01}=2,\qquad T_{10}=3.
\]
Note that the product of noncommuting elements cannot be equal to one of them. Thus, if for a given order there is only one nonabelian group (this is true for orders \(6,10,14,21,22\)), then two filled cells are sufficient for the unique identification of the group.
\end{proof}

The property ``to be cyclic'' is expressed in the most universal way from the point of view of the order of the group: in this case the group contains an element whose order is not a proper divisor of \(n\).

Let us describe a universal algorithm for constructing a record of this assertion with the help of Cayley tables.

\begin{prop}
For the assertion ``the group has an element of order at least \(n\)'', it is sufficient to fill \(n\) cells of the Cayley table.
\end{prop}

\begin{proof}
Denote the element with number \(0\) by \(a\). Fill the Cayley table as follows:
\[
T_{10}=1,\quad T_{01}=2,\quad\ldots,\quad T_{0,n-1}=n.
\]
The first entry is equivalent to the assertion \(a\alpha=a^2\ne a\) (from this it follows that \(a\ne e\)). From the second entry it follows that \(a\alpha^2=a^3\ne a\), i.e.\ \(\operatorname{ord}(a)\ne2\). Continuing the arguments in a similar way, we obtain that the element \(a\) has order at least \(n\).
\end{proof}

The proved assertion gives an upper bound for the number of filled cells of the Cayley table necessary for the unique identification of a group of order \(n\). Since the possible orders of an element of a group are divisors of the order of the group, not all equalities written in Proposition 2 are necessary for unique identification.

\begin{prop}
The assertion ``the group of order \(2^k\) is cyclic'' can be written by filling \(k\) cells of the Cayley table.
\end{prop}

\begin{proof}
Denote a generator by \(a\). Start the enumeration of the elements of the group with
\[
a,a^2,a^4,\ldots,a^{2^{k-1}}.
\]
Fill the Cayley table as follows:
\[
T_{00}=1,\quad T_{11}=2,\quad\ldots,\quad T_{k-1,k-1}=k.
\]
Thus,
\[
a^2a^2=a^4\ne a^2\quad(\operatorname{ord}(a)\ne2),
\]
\[
\ldots,
\]
\[
a^{2^{k-1}}a^{2^{k-1}}=a^{2^k}\ne a^{2^{k-1}}
\quad(\operatorname{ord}(a)\ne2^{k-1}).
\]
Therefore, the order of the element \(a\) is equal to \(2^k\), but an element of such order can exist only in the cyclic group \(\mathbb{Z}_{2^k}\). Consequently, the indicated restrictions uniquely determine the group.
\end{proof}

In the general case, the idea of constructing a partial filling remains the same: we enumerate the possible proper divisors of the order of the group and choose a filling showing that the given power of the generator is a new element of the group.

\begin{prop}
The cyclic group of order \(9\) can be uniquely identified by two filled cells of the Cayley table, and the group of order \(12\) by three cells.
\end{prop}

\begin{proof}
In the case of the group \(\mathbb{Z}_9\), start the enumeration of the elements of the group with \(a,a^2,a^4\) (where \(a\) is a generator of \(\mathbb{Z}_9\)). Filling
\[
T_{00}=1,\qquad T_{11}=2
\]
shows that
\[
aa=a^2\ne a
\quad\text{and}\quad
a^2a^2=a^4\ne a.
\]

In the case of the group \(\mathbb{Z}_{12}\), start the enumeration of the elements of the group with \(a,a^2,a^4,a^8\) (where \(a\) is a generator of \(\mathbb{Z}_{12}\)). Filling
\[
T_{00}=1,\qquad T_{11}=2,\qquad T_{22}=3
\]
shows that
\[
aa=a^2\ne a,\qquad
a^2a^2=a^4\ne a^2,a
\quad(\operatorname{ord}(a)\ne2,3),
\]
and
\[
a^4a^4=a^8\ne a^2,a^4
\quad(\operatorname{ord}(a)\ne4,6).
\]
\end{proof}

In some cases it is possible to construct a more economical partial filling for noncyclic groups. In particular, such fillings exist for groups of order \(18\) and \(20\).

\begin{prop}
For the group \(S_3\times\mathbb{Z}_3\) (a group of order \(18\)) there exists a filling of \(4\) cells of the Cayley table which uniquely identifies it.
\end{prop}

\begin{proof}
In the remaining groups of order \(18\), a nonidentity element can be a square of at most two elements of the group, whereas in the group \(S_3\times\mathbb{Z}_3\) the following is true:
\[
\langle(12),a\rangle\langle(12),a\rangle
=\langle e,a^2\rangle\ne\langle e,e\rangle,
\]
\[
\langle(13),a\rangle\langle(13),a\rangle
=\langle e,a^2\rangle\ne\langle e,e\rangle,
\]
\[
\langle(23),a\rangle\langle(23),a\rangle
=\langle e,a^2\rangle\ne\langle e,e\rangle,
\]
where \(\langle f,g\rangle\) is the ordered pair \(f\in S_3\), \(g\in\mathbb{Z}_3\); \(e\) is the identity of the groups \(S_3\) and \(\mathbb{Z}_3\); \((12)\in S_3\) is a transposition; and \(a\) is a generator of \(\mathbb{Z}_3\).
\end{proof}

Among the groups of order \(20\), the minimal number of cells is required for the unique identification of the group \(GA(1,5)\) (the group of one-dimensional affine transformations over the field of \(5\) elements).

\begin{prop}
For the group \(GA(1,5)\) (a group of order \(20\)) there exists a filling of \(3\) cells of the Cayley table which uniquely identifies it.
\end{prop}

\begin{proof}
The group \(GA(1,5)\) is the unique nonabelian group of order \(20\) with elements satisfying the relation
\[
abab=ba.
\]
Accordingly, we can start the enumeration as follows:
\[
a,b,ab,ba,\ldots
\]
and fill three cells as
\[
a\cdot b=ab,\qquad b\cdot a=ba,\qquad ab\cdot ab=ba.
\]
\end{proof}

\section{Checking minimality}

It remains only to check that there are no uniquely identifying restrictions with a smaller number of cells.

These calculations were performed with the help of the computer algebra system GAP [5], using the module SmallGrp [6].

\begin{enumerate}
\item From the corresponding library, the list of Cayley tables of groups of the required order is imported up to isomorphism.
\item A variant of a partial filling of the Cayley table (called below a constraint) is chosen.
\item All possible permutations are applied to this constraint (since this is the same constraint but with another numbering of the elements of the group).
\item By brute force over the Cayley tables, it is established how many of them satisfy the set of constraints (i.e.\ under some renumbering of the elements of the group the original constraint will hold).
\end{enumerate}

The calculations are presented in the repository [7].

\section{Results for groups of orders \(2\)--\(23\)}

Thus we obtain Table~\ref{tab:main}, which gives the minimal numbers of filled cells of the Cayley table for the unique identification of groups of order \(n\le 23\).

\begin{table}[htbp]
\centering
\caption{Minimal number of Cayley table cells necessary for the unique identification of a group}
\label{tab:main}
\begin{tabular}{|c|c|}
\hline
Order of the group & Minimum number of filled cells \\
\hline
\(1,2,3,5,7,11,13,15,17,19,23\) & \(0\) \\
\(4,6,9,10,14,21,22\) & \(2\) \\
\(8,12,20\) & \(3\) \\
\(16,18\) & \(4\) \\
\hline
\end{tabular}
\end{table}

\section{Conclusion}

In this article we have studied the question of the minimal number of cells of the Cayley table whose filling allows one to uniquely identify a group. Minimal partial fillings have been constructed for groups of order \(\le 23\).

As a result of the research, an integer sequence was obtained and included in the Online Encyclopedia of Integer Sequences (OEIS) under the number A387855 [8].

Further research can develop in the following directions: computation of the minimal filling for higher orders, computation of a universal filling allowing one to identify all groups of a given order. It is also possible to study a more general quantity, namely the number of nonisomorphic groups of order \(n\) when \(k\) cells of the Cayley table are filled, which represents a so-called ``triangle of sequences'' (on triangles of sequences see [9]).

\end{document}